\documentclass[12pt,reqno]{amsart}
\usepackage{amsmath,amssymb,amsfonts,amscd,latexsym,amsthm,mathrsfs}
\usepackage[unicode]{hyperref}
\renewcommand{\d}{{\mathrm d}}

\newcommand{\mm}{{\mathrm m}}

\begin{document}

\title{Mahler measure numerology}

\author{Wadim Zudilin}
\address{Department of Mathematics, IMAPP, Radboud University, PO Box 9010, 6500\,GL Nijmegen, Netherlands}
\email{w.zudilin@math.ru.nl}

\date{17 September 2021}

\dedicatory{To David Boyd, an oustanding conductor of the Mahler measures,\\on the ocassion of his 80th birthday}

\begin{abstract}
We discuss some (conjectural) evaluations of $L$-values attached to elliptic curves of conductors 15, 21, 24 and 32
as `hypergeometric periods'.
These numerical observations are motivated by the Mahler measures of three-variable polynomials.
\end{abstract}


\maketitle
%==================================================

The (logarithmic) Mahler measure
\begin{align*}
\mm(P(x_1,\dots,x_k))
&=\frac1{(2\pi i)^k}\idotsint_{|x_1|=\dots=|x_k|=1}\log|P(x_1,\dots,x_k)|\,\frac{\d x_1}{x_1}\dotsb\frac{\d x_k}{x_k}
\end{align*}
of an $k$-variable (Laurent) polynomial $P(x_1,\dots,x_k)\in\mathbb C[x_1,\dots,x_k]$
is a quite unique attractor of numerous problems in mathematics \cite{BZ20}.
One of big problems (well open even in the case $k=1$!) is the range of the Mahler measures attached to polynomials with \emph{integer} coefficients \cite{Bo81}.
A particular aspect of this problem is a remarkable connection of such Mahler measures to the $L$-values of algebraic varieties, usually related to the zero loci of the underlying polynomials $P(x_1,\dots,x_k)$.
Research in this direction originated in the 1990s: the works of Deninger \cite{De97}, Boyd \cite{Bo98} and Rodriguez Villegas \cite{Vi99} set up in depth the story for the case of two-variable polynomials $P(x,y)$ such that $P(x,y)=0$ defines an elliptic curve.
The papers all feature the family of two-variable Mahler measures
$$
\mu(k)=\mm\bigl(k+(x+1/x)(y+1/y)\bigr),
$$
which is indeed quite special.
Illustrative examples for the range $0<k<4$ include
\begin{alignat*}{3}
\mu(1)&=L'(E_{15a8},0), &\quad
\mu(\sqrt2)&=4L'(E_{56a1},0), &\quad
\mu(2)&=L'(E_{24a4},0), \\
\mu(2\sqrt2)&=L'(E_{32a1},0), &\quad
\mu(3)&=\tfrac12L'(E_{21a4},0), &\quad
\mu(2\sqrt3)&\overset?=2L'(E_{72a1},0),
\end{alignat*}
where $E_{NxM}$ refers to the elliptic curve labeled in accordance with \cite{LMFDB}
(the label corresponds to a Weierstrass form of the curve $k+(x+1/x)(y+1/y)=0$ which is defined over $\mathbb Q$ whenever $k^2\in\mathbb Q\setminus\{0,16\}$).
Notice that for this range
\begin{equation}
\mu(k)
=\frac k4\cdot{}_3F_2\biggl(\begin{matrix} \tfrac12, \, \tfrac12, \, \tfrac12 \\ 1, \, \tfrac32 \end{matrix}\biggm| \frac{k^2}{16} \biggr)
=\int_0^{k/4}{}_2F_1\biggl(\begin{matrix} \tfrac12, \, \tfrac12 \\ 1 \end{matrix}\biggm| x^2 \biggr)\,\d x,
\label{mu2}
\end{equation}
so that the Mahler measure evaluations can be translated into equalities of the $L$-values and the values of hypergeometric functions.

The $L$-values of conductor 15, 21 and 24 elliptic curves also show up, conjecturally(!), in the three-variable Mahler measures%
\footnote{The notation `$\overset?=$' is used for equalities numerically observed.}
\begin{align*}
\mm(1+x+y+xy+z)&\overset?=-2L'(E_{15},-1)=\frac{15^2}{4\pi^4}L(E_{15},3)=0.4839979734\hdots,
\\
\mm((1+x)^2+y+z)&\overset?=-L'(E_{24},-1)=\frac{72}{\pi^4}L(E_{24},3)=0.7025655062\hdots,
\\
\mm(1+x+y-xy+z)&\overset?=-\tfrac54L'(E_{21},-1)=\frac{5\cdot21^2}{32\pi^4}L(E_{21},3)=0.6242499823\hdots\,.
\end{align*}
The Mahler measures in all these cases are reduced to hypergeometric integrals \cite[Section~6.3]{BZ20} (see also \cite{SZ20}), and there are also hypergeometric expressions available for $L'(E_{32},-1)$ \cite{Zu13}.
Here we highlight the corresponding formulae for the $L$-values, dropping off the intermediate Mahler measures, and give some additional formulae that are not linked to known conjectures on the Mahler measures but reflect similarities with the two-variable expression~\eqref{mu2}.

\medskip
\noindent
\emph{Conductor $15$}  ($k=1$):
The expectation in this case is
\begin{align}
-L'(E_{15},-1)
&\overset?=\frac1{\pi}\int_{1/4}^1{}_2F_1\biggl(\begin{matrix} \tfrac12, \, \tfrac12 \\ 1 \end{matrix}\biggm| 1-x^2 \biggr)\log(4x)\,\d x
\nonumber\\
&=-\frac1{\pi}\int_0^{1/4}{}_2F_1\biggl(\begin{matrix} \tfrac12, \, \tfrac12 \\ 1 \end{matrix}\biggm| 1-x^2 \biggr)\log(4x)\,\d x.
\label{eq:L15-1}
\end{align}
Notice the induced evaluation
\begin{equation}
\frac1{\pi}\int_0^1{}_2F_1\biggl(\begin{matrix} \tfrac12, \, \tfrac12 \\ 1 \end{matrix}\biggm| 1-x^2 \biggr)\log(4x)\,\d x
=0
\label{eq:L15-1a}
\end{equation}
established in~\cite{SZ20}.
There is a complimentary observation of H.~Cohen (2018) in this case, not related to a Mahler measure:
\begin{equation}
L'(E_{15},-2)
\overset?=\frac{12}{\pi^2}\int_0^{1/4}{}_2F_1\biggl(\begin{matrix} \tfrac12, \, \tfrac12 \\ 1 \end{matrix}\biggm| x^2 \biggr)\log^2(4x)\,\d x.
\label{eq:L15-2}
\end{equation}

\medskip
\noindent
\emph{Conductor $21$}  ($k=3$):
Here we expect
\begin{equation}
-5L'(E_{21},-1)
\overset?=\frac2{\pi}\int_{1/4}^2{}_2F_1\biggl(\begin{matrix} \tfrac12, \, \tfrac12 \\ 1 \end{matrix}
\biggm| 1-(1-x)^2 \biggr)\log(4x)\,\d x.
\label{eq:L21-1}
\end{equation}
The argument $2x-x^2$ of the hypergeometric integrand in \eqref{eq:L21-1} ranges between $1-k^2/16$ and $1$ when $\frac14<x<1$ and then between $1$ and $0$ when $1<x<2$. Thus, taking into account the evaluation%
\footnote{This is kindly reported to us by Ringeling~\cite{Ri20b} based on the findings in his work~\cite{Ri20a}.}
\begin{equation}
\frac{7\zeta(3)}{\pi^2}=-28\zeta'(-2)
=\frac1{\pi}\int_0^1{}_2F_1\biggl(\begin{matrix} \tfrac12, \, \tfrac12 \\ 1 \end{matrix}
\biggm| 1-x^2 \biggr)\log(4(1+x))\,\d x
\label{eq:L21-1a}
\end{equation}
for the second part (see also \eqref{eq:L15-1a}), our main expectation is
\begin{equation}
-\frac52L'(E_{21},-1)+28\zeta'(-2)
\overset?=\frac1{\pi}\int_0^{3/4}{}_2F_1\biggl(\begin{matrix} \tfrac12, \, \tfrac12 \\ 1 \end{matrix}
\biggm| 1-x^2 \biggr)\log(4(1-x))\,\d x.
\label{eq:L21-1b}
\end{equation}

\medskip
\noindent
\emph{Conductor $24$} ($k=2$):
We know in this case \cite{SZ20} that
\begin{align}
-L'(E_{24},-1)
&\overset?=\frac{2\sqrt2}{\pi\sqrt\pi}
\sum_{n=0}^\infty\frac{(\frac12)_n\,\Gamma(\frac n2+\frac34)}{n!\,\Gamma(\frac n2+\frac54)\,(2n+1)^2}\,\Bigl(\frac12\Bigr)^n
\label{eq:L24-1}
\\
&=\frac{8\Gamma(\frac34)^2}{\pi^{5/2}}\,
{}_5F_4\biggl(\begin{matrix} \tfrac14, \, \tfrac14, \, \tfrac14, \, \tfrac34, \, \tfrac34 \\ \frac12, \, \tfrac54, \, \tfrac54, \, \tfrac54 \end{matrix}\biggm| \frac 14 \biggr)
+\frac{\Gamma(\frac14)^2}{54\pi^{5/2}}\,
{}_5F_4\biggl(\begin{matrix} \tfrac34, \, \tfrac34, \, \tfrac34, \, \tfrac54, \, \tfrac54 \\ \frac32, \, \tfrac74, \, \tfrac74, \, \tfrac74 \end{matrix}\biggm| \frac 14 \biggr).
\nonumber
\end{align}
However there is an identification of $L'(E_{24},-1)$ similar to \eqref{eq:L15-1} and \eqref{eq:L21-1b} but not related to a Mahler measure:
\begin{equation}
-2L'(E_{24},-1)+28\zeta'(-2)
\overset?=\frac1{\pi}\int_0^{2/4}{}_2F_1\biggl(\begin{matrix} \tfrac12, \, \tfrac12 \\ 1 \end{matrix}
\biggm| 1-x^2 \biggr)\log\Bigl(\frac{1-x}x\Bigr)\,\d x.
\label{eq:L24-1a}
\end{equation}
In a different direction, here is a striking connection of \eqref{eq:L24-1} with other $L$-values of the same elliptic curve%
\footnote{The proof of \eqref{eq:L24-a} we are aware of is somewhat tricky. One starts with
\begin{align*}
&
\frac{\sqrt\pi}2\sum_{n=0}^\infty\frac{(\frac12)_n\,\Gamma(\frac n2+\frac34)}{n!\,\Gamma(\frac n2+\frac54)}\,z^n
=\sum_{n=0}^\infty\frac{(\frac12)_n}{n!}\,z^n\int_0^1t^{n+1/2}(1-t^2)^{-1/2}\,\d t
%\\ &\quad
=\int_0^1\frac{\sqrt t}{\sqrt{(1-t^2)(1-zt)}}\,\d t
\\ &\quad
=\sqrt2\,\bigg(\int_0^1\frac{\d u}{(1-\frac12u^2)\sqrt{(1-u^2)(1-\frac12(1+z)u^2)}}
-\int_0^1\frac{\d u}{\sqrt{(1-u^2)(1-\frac12(1+z)u^2)}}\bigg)
\end{align*}
to arrive, when $z=\frac12$, at
$$
\int_0^1\frac{\d u}{(1-\frac12u^2)\sqrt{(1-u^2)(1-\frac34u^2)}}
=\frac32\int_0^1\frac{\d u}{\sqrt{(1-u^2)(1-\frac34u^2)}},
$$
and this allows one to conclude with \eqref{eq:L24-a} from \eqref{eq:L24'-a}.
There seems to be no simple relation like this for $z=-\frac12$ (and for other values of $z$ from the real interval $|z|<1$). This fact makes it difficult to lift the finding to \eqref{eq:L24-0}.}:
\begin{align}
L'(E_{24},0)
&\overset?=\frac1{\sqrt{2\pi}}
\sum_{n=0}^\infty\frac{(\frac12)_n\,\Gamma(\frac n2+\frac34)}{n!\,\Gamma(\frac n2+\frac54)\,(2n+1)}\,\Bigl(-\frac12\Bigr)^n,
\label{eq:L24-0}
\\
L(E_{24},1)
&=\frac{\sqrt\pi}{4\sqrt2}
\sum_{n=0}^\infty\frac{(\frac12)_n\,\Gamma(\frac n2+\frac34)}{n!\,\Gamma(\frac n2+\frac54)}\,\Bigl(\frac12\Bigr)^n.
\label{eq:L24-a}
\end{align}
These two evaluations can be compared with
\begin{align}
L'(E_{24},0)
&=\frac12\sum_{n=0}^\infty\frac{(\frac12)_n^2}{n!^2\,(2n+1)}\,\Bigl(\frac14\Bigr)^n,
\label{eq:L24'-0}
\\
L(E_{24},1)
&=\frac\pi8\sum_{n=0}^\infty\frac{(\frac12)_n^2}{n!^2}\,\Bigl(1-\frac14\Bigr)^n,
\label{eq:L24'-a}
\end{align}
where the first one comes out of the formula for $\mu(2)$ interpreted through \eqref{mu2}.

\medskip
\noindent
\emph{Conductor $32$} ($k=2\sqrt2$):
Here we have the following \cite{Zu13}:
\begin{align}
-L'(E_{32},-1)
&=\frac{\Gamma(\frac14)^2}{6\sqrt2\pi^{5/2}}\,
{}_4F_3\biggl(\begin{matrix} 1, \, 1, \, 1, \, \frac12 \\ \frac74, \, \frac32, \, \frac32 \end{matrix}\biggm| 1 \biggr)
+\frac{4\Gamma(\frac34)^2}{\sqrt2\pi^{5/2}}\,
{}_4F_3\biggl(\begin{matrix} 1, \, 1, \, 1, \, \frac12 \\ \frac54, \, \frac32, \, \frac32 \end{matrix}\biggm| 1 \biggr)
\nonumber\\ &\qquad
+\frac{\Gamma(\frac14)^2}{2\sqrt2\pi^{5/2}}\,
{}_4F_3\biggl(\begin{matrix} 1, \, 1, \, 1, \, \frac12 \\ \frac34, \, \frac32, \, \frac32 \end{matrix}\biggm| 1 \biggr).
\label{eq:L32-2}
\end{align}
Conjecturally we find out the following representation resembling \eqref{eq:L15-1}, \eqref{eq:L21-1b} and \eqref{eq:L24-1a}:
\begin{equation}
-L'(E_{32},-1)
\overset?=\frac1{\pi}\int_0^{2\sqrt{2}/4}{}_2F_1\biggl(\begin{matrix} \tfrac12, \, \tfrac12 \\ 1 \end{matrix}
\biggm| 1-x^2 \biggr)\log\biggl(\frac{1-x^2}{x^2}\biggr)\,\d x.
\label{eq:L32-1a}
\end{equation}

\medskip
\noindent
\textbf{Acknowledgements.}
I thank Fran\c cois Brunault and Berend Ringeling for inspiring conversations and assistance.

%==================================================

\end{document}